\title{
Minimizing the number of carries in addition
}
\author{Noga Alon
\thanks{Sackler School of Mathematics
and Blavatnik School of Computer Science,
Tel Aviv University,
Tel Aviv 69978, Israel.
Email: {\tt nogaa@tau.ac.il}.
Research supported in part by an ERC Advanced
grant, by a USA-Israeli BSF grant and by 
the Israeli I-Core program.
}
}

\documentclass[11pt]{article}
\usepackage{amsfonts,amssymb,amsmath,latexsym}
\def\qed{\ifvmode\mbox{ }\else\unskip\fi\hskip 1em plus 10fill$\Box$}

\newtheorem{theo}{Theorem}[section]

\newtheorem{coro}[theo]{Corollary}

\begin{document}
\date{}

\maketitle

\begin{abstract}
When numbers are added in base $b$ in the usual way, carries occur.
If two random, independent 1-digit numbers are added, then the
probability of a carry is $\frac{b-1}{2b}$. Other choices of digits lead
to less carries. In particular, if for odd $b$ we  use the  digits
$\{-(b-1)/2, -(b-3)/2, \ldots , \ldots (b-1)/2\}$
then the
probability of carry is only $\frac{b^2-1}{4b^2}$.
Diaconis, Shao and Soundararajan conjectured that this 
is the best choice of digits, and proved that this is
asymptotically the case when $b=p$ is a large prime. In this note 
we prove this conjecture for all odd primes $p$.
\vspace{0.1cm}

\noindent
AMS Subject Classification: 11P99
\vspace{0.1cm}

\noindent
Keywords: carry, modular addition

\end{abstract}
\section{The problem and result}

When numbers are added in base $b$ in the usual way, carries occur.
If two added one-digit numbers are random  and independent, then the
probability of a carry is $\frac{b-1}{2b}$. Other choices of digits lead
to less carries. In particular, if for odd $b$ we  use 
{\em balanced digits}, that is, the  digits
$$\{-(b-1)/2, -(b-3)/2, \ldots ,0, 1, 2, \ldots (b-1)/2\}$$
then the
probability of carry is only $\frac{b^2-1}{4b^2}$.
Diaconis, Shao and Soundararajan \cite{DSS} conjectured that this 
is the best choice of digits, and proved that this conjecture is
asymptotically correct when $b=p$ is a large prime. More precisely,
they proved the following.
\begin{theo}[\cite{DSS}]
\label{t11}
For every $\epsilon >0$ there exists a  number $p_0=p_0(\epsilon)$
so that for any prime $p>p_0$ the probability of carry when 
adding two random independent one-digit numbers using
any fixed set of digits in base $p$ is at least
$\frac{1}{4}-\epsilon$.
\end{theo}
The estimate given in \cite{DSS} to $p_0=p_0(\epsilon) $ is a tower
function of $1/\epsilon$.

Here we establish a tight result for any prime $p$, proving the
conjecture for any prime.
\begin{theo}
\label{t12}
For any odd prime $p$, the probability of carry when 
adding two random independent one-digit numbers using
any fixed set of digits in base $p$ is at least
$\frac{p^2-1}{4p^2}$.
\end{theo}
The proof is very short, and is in fact mostly an observation that
the above result follows from a theorem of  J. M. Pollard proved
in the 70s. The conjecture for non-prime values of $p$ remains open.

The proof can be extended to show that
balanced digits minimize the probability 
of carry while adding $k$ numbers, for any $k \geq 2$.

The rest of this short note is organized as follows. The next
section  contains a brief description of the digit systems
considered. In Section 3 we present the short derivation of Theorem 
\ref{t12} from Pollard's Theorem, and in Section 4 we describe the
proof of the extension to the addition of more than two summands.

\section{Addition and choices of digits}

A simple example illustrating the advantage of using digits that
minimize the probability of carry is that of adding numbers in the
finite cyclic group $Z_{b^2}$. Here the basis used is $b$. 
Since $Z_{b^2}$ is a finite
group, one can choose random members of it $g_1,g_2, \ldots ,g_n$ 
uniformly and independently and consider their sum in $Z_{b^2}$. 
Following
the discussion in Section 6 of \cite{BDF}, consider the normal
subgroup $ Z_b \triangleleft Z_{b^2},$
where $Z_{b}$  is the subgroup of $Z_{b^2}$ consisting of the
elements $\{0,b, 2b, \ldots ,(b-1)b\}$.

Let $A \subset Z_{b^2}$ be a set of representatives of the 
cosets of $Z_b$ in $Z_{b^2}$. Therefore $|A|=b$ and no two elements
of $A$ are equal modulo $b$. These are the digits we use.
Any element $g \in Z_{b^2}$ now has a unique representation of the
form $g=x+y$, where $x \in A$ and $y \in Z_b$. Indeed, $x$ is the
member of $A$ representing the coset of $Z_b$ that contains 
$g$, and $y \in Z_b$ is determined by the equality $g=x+y$.

Suppose, now that $g_i=x_i+y_i$ with $x_i \in A$ and $y_i \in Z_b$
is the representation of $n$ elements $g_i$ of $Z_{b^2}$ that we
wish to sum. We start by computing $g_1+g_2$. To do so one first
adds $x_1+x_2$ (in $Z_{b^2}$). 
If their sum, call it $z_2$, is a member of $A$,
then there is no carry in this stage. We can then compute the sum 
$y_1+y_2$ in $Z_b$, and get $w_2$ (in this stage there is no
carry, as the addition here is modulo $b$). Therefore, in this case
the representation of $g_1+g_2$ using our digits is
$g_1+g_2=z_2+w_2$ with $z_2 \in A$ and $w_2 \in Z_p$, and we can
now proceed by induction and compute the sum of this element
with $g_3+g_4 + \ldots +g_n$.  If, on the other hand,
$x_1+x_2 \not \in A$ then there is a carry. In this case we let
$z_2$ be the unique member of $A$ so that $x_1+x_2$ lies in the
coset of $Z_b$ containing $z_2$. This determines the element
$t_2 \in Z_b$ so that $x_1+x_2=z_2+t_2$. The carry here is $t_2$,
and we can now proceed and compute the sum $t_2+y_1+y_2$ in
$Z_b$ getting an element $w_2$. Therefore in this case too
the unique representation of $g_1+g_2$ using our digits is
$z_2+w_2$, but since the process of computing them involved the
carry $t_2$ the number of additions performed during the
computation in the group $Z_b$ was $2$ and not $1$ as in the case
that involved no carry.

As $g_1$ and $g_2$ are random independent elements chosen uniformly
in $Z_{b^2}$, their sum is also uniform in this group, implying
that the element $z_2+w_2$ is also uniform. Therefore, when we now
proceed and compute the sum of this element with $g_3$ the
probability of carry is again exactly as it has been before (though
the conditional probability of a carry given that there has been
one in the previous step  may be different). We conclude that the
expected number of carries during the whole process of adding
$g_1+g_2+ \cdots +g_n$ that consist of $n-1$ additions is exactly
$(n-1)$ times the probability of getting a carry in one addition of
two independent uniform random digits in the set $A$.

Therefore, the problem of minimizing the expected number of these carries 
is that of
selecting the set of coset representatives $A$ 
so that the probability that the addition of two random members of
$A$ in $Z_{b^{2}}$  does not lie  in $A$ is minimized.

This leads to the following equivalent
formulation of Theorem \ref{t12}.
\begin{theo}
\label{t21}
Let $p$ be an odd prime.
For any subset $A$ of the group $Z_{p^2}$ of integers modulo $p^2$
so that $|A|=p$ and the members of $A$ are pairwise distinct modulo
$p$, the number of ordered pairs $(a,b) \in A \times A$ so that
$a+b~(mod~p^2) \not \in A$ is at least $\frac{p^2-1}{4}$.
\end{theo}

\section{Adding two numbers}

The result of Pollard needed here is the following. 
\begin{theo}[Pollard \cite{Po}]
\label{t31}
For an integer
$m$ and two sets $A$ and $B$ of residues modulo $m$, and for any
positive integer $r$, let $N_r=N_r(A,B)$ denote the number of all
residues modulo $m$ that have a representation as a sum
$a+b$ with $a \in A$ and $b \in B$ in at least $r$ ways.
(The representations are counted as ordered pairs, that is, if
$a$ and $b$ differ and both belong to $A \cap B$, then 
$a+b=b+a$ are two distinct representations of the sum).
If $(x-y,m)=1$ for any two distinct elements $x,y \in B$ then
for any $1 \leq r \leq \min\{|A|,|B|\}$:
$$
N_1+N_2 + \cdots +N_r \geq r \cdot \min\{m,|A|+|B|-r\}.
$$
\end{theo}
\vspace{0.3cm}

\noindent
Note that the case $r=1$ is the classical theorem of Cauchy and
Davenport (see \cite{Ca}, \cite{Da}). The proof is short and
clever, following the approach in the original papers of Cauchy and
Davenport. It proceeds by induction on $|B|$, where in the induction
step one first replaces $B$ by a shifted  copy $B'=B-g$ so that
$I=A \cap B'$ satisfies $1 \leq |I| < |B'|=|B|$, and then applies
the induction hypothesis to the pair $(A \cup B', I)$ and to the pair
$(A-I, B'-I)$. The details can be found in \cite{Po} (see also
\cite{Po0}). 
\vspace{0.3cm}

\noindent
{\bf Proof of Theorem \ref{t12}:}\,
As mentioned above, the statement of Theorem \ref{t12}
is equivalent to  that of Theorem \ref{t21}:
for any subset $A$ of the group $Z_{p^2}$ of integers modulo $p^2$
so that $|A|=p$ and the members of $A$ are pairwise distinct modulo
$p$, the number of ordered pairs $(a,b) \in A \times A$ so that
$a+b~(mod~p^2) \not \in A$ is at least $\frac{p^2-1}{4}$. 

Given such a set $A$, note that $(x-y,p^2)=1$ for every two
distinct elements $x,y \in A$. We can thus apply Pollard's Theorem
stated above
with $m=p^2$, $A=B$ and $r=(p-1)/2$ to conclude that 
$$
N_1+N_2+ \cdots +N_r \geq r \cdot \min\{p^2,|A|+|B|-r\}=
r \cdot (2p-r).
$$
The sum $N_1+N_2 + \cdots + N_r$ counts every element $x \in Z_{p^2}$
exactly $\min\{r,n(x)\}$ times, where $n(x)$ is the number of
representations of $x$ as an ordered sum $a+b$ with $a,b \in A$.
The total contribution to this sum arising from elements $x \in A$
is at most $r |A|=rp$. Therefore, there are at least
$r(p-r)=\frac{b-1}{2}\frac{b+1}{2}=\frac{b^2-1}{4}$ ordered pairs
$(a,b) \in A \times A$ so that $a+b \not \in A$, completing the
proof.  $\Box$
\vspace{0.3cm}

\noindent
{\bf Remark:}
The above proof shows that even if we use one set of digits for one
summand, a possibly different set of digits for the second summand,
and a third set of digits for the sum, then still the probability of
carry must be at least $\frac{b^2-1}{4b^2}$.

\section{Adding more numbers}

The theorem of Pollard is more general than the
statement above and deals with addition of $k>2$ sets as well. This
can be used in
determining the minimum possible probability of carry
in the addition  of $k$ random 1-digit numbers in a prime base $p$ with
the best choice of the  $p$ digits. In fact, for every $k$ 
and every odd prime $p$, the
minimum probability is obtained by using the balanced digits
$\{-(p-1)/2, -(p-3/2, \ldots ,(p-1)/2\}$. 
Therefore, the minimum possible probability of
carry in adding $k$ 1-digit numbers in a prime  base $p>2$  is exactly the
probability that the sum of $k$ independent random variables, each
distributed uniformly on the set
$$
\{-(p-1)/2, -(p-3/2, \ldots ,(p-1)/2\},
$$
is of absolute value exceeding $(p-1)/2$.

Here are the details. We need the following.
\begin{theo}[Pollard \cite{Po}]
\label{t41}
Let $m$ be a positive integer and let $A_1, A_2, \ldots ,A_k$ be
subsets of $Z_m$. Assume, further, that for every $2 \leq i \leq k$
every two distinct elements $x,y$ of $A_i$ satisfy $(x-y,m)=1$. 
Let $A'_1,A'_2, \ldots ,A'_k$ be another collection of subsets of
$Z_m$, in which each $A'_i$ consists of consecutive elements
and satisfies $|A_i'|=|A_i|$. For an $x \in Z_m$ let $n(x)$ denote
the number of representations of $x$ as an ordered sum of the form
$x=a_1+a_2 + \ldots +a_k$ with $a_i \in A_i$, and let
$n'(x)$ denote the number of representations of $x$ as an ordered
sum of the form $x =a'_1+a'_2 + \cdots +a'_k$ with $a_i' \in A'_i$ 
for all $i$. Then, for any integer $r \geq 1$
$$
\sum_{x \in Z_m} \min\{r,n(x)\} \geq \sum_{x \in Z_m} \min
\{r,n'(x) \}.
$$
\end{theo}
\begin{coro}
\label{c42}
Let $p$ be an odd prime, and let $A$ be a subset of cardinality $p$
of $Z_{p^2}$. Assume,
further, that the members of $A$ are pairwise distinct modulo $p$.
Put $A'=\{-(p-1)/2, -(p-3)/2, \ldots ,(p-1)/2\}$. Then, for any
positive integer $k$, the
number of ordered sums modulo $p^2$ of $k$ elements of $A$ that do
not belong to $A$ is at least as large as the 
number of ordered sums modulo $p^2$ of $k$ elements of $A'$ that do
not belong to $A'$.
\end{coro}
\vspace{0.3cm}

\noindent
{\bf Proof:}\, Let $r$ be the number of  ordered sums modulo $p^2$
of $k$ elements of $A'$ whose value is precisely $(p-1)/2$. It is
not difficult to check that for any other member $g$ of $A'$ there
are at least $r$ ordered sums modulo $p^2$
of $k$ elements of $A'$ whose value is precisely $g$. 
Similarly, for any $x  \not \in A'$, 
the number of ordered sums of $k$ elements of $A'$ whose value is
precisely $x$ is at most $r$. Indeed, the number of times an
element is obtained as an ordered sum modulo $p^2$ 
of $k$ elements of $A'$ 
is a monotone non-increasing function of its distance from $0$.
(This  can be easily proved by induction on $k$).
Therefore, 
$$
\sum_{x \in Z_{p^2}} \min \{r,n'(x)\} =r p +
\sum_{x \in Z_{p^2}-A'} n'(x).
$$
By Theorem \ref{t41} with $m=p^2$, $A_1=A_2= \ldots =A_k=A$ and the value of
$r$ above
$$
\sum_{x \in Z_{p^2}} \min \{r,n(x)\} \geq
\sum_{x \in Z_{p^2}} \min \{r,n'(x)\}.
$$
However, clearly,
$$
rp + \sum_{x \in Z_{p^2}-A} n(x)
\geq  \sum_{x \in Z_{p^2}} \min \{r,n(x)\},
$$
and therefore
$$
rp + \sum_{x \in Z_{p^2}-A} n(x)
\geq 
\sum_{x \in Z_{p^2}} \min \{r,n'(x)\}=
r p +
\sum_{x \in Z_{p^2}-A'} n'(x).
$$
This implies that
$$
\sum_{x \in Z_{p^2}-A} n(x)
\geq 
\sum_{x \in Z_{p^2}-A'} n'(x),
$$
as needed.  $\Box$

The corollary clearly implies that
the minimum possible probability of
carry in adding $k$ $1$-digit numbers in a prime  base $p>2$  is exactly the
probability that the sum of $k$ independent random variables, each
distributed uniformly on the set
$$
\{-(p-1)/2, -(p-3/2, \ldots ,(p-1)/2\},
$$
is of absolute value exceeding $(p-1)/2$.
\vspace{0.2cm}

\noindent
{\bf Remarks:}\,

\begin{itemize}
\item
As in the case of two summands, the proof implies that
the assertion of the last paragraph holds even if we
are allowed to choose a different set of digits for each summand
and for the result. 
\item
After the completion of this note I learned from the authors of
\cite{DSS} that closely related results (for addition in $Z_p$,
not in $Z_{p^2}$) appear in the paper of Lev \cite{Le}.
\end{itemize}
\vspace{0.3cm}

\noindent
{\bf Acknowledgment}\,
Part of this work was carried out during a visit at 
the Georgia Institute of Technology.
I would like to thank my hosts
for their hospitality during this visit. 
I would also like to thank Persi Diaconis for an inspiring lecture
in Georgia Tech in which he mentioned the problem considered here
and Kannan Soundararajan for telling me about \cite{Le}.

\end{document}